\documentclass[11pt]{article}
\usepackage{epsfig,makeidx,amsmath,amsfonts,latexsym,subfigure}
\newtheorem{theorem}{Theorem}[section]
\newtheorem{prop}{Proposition}[section]
\newtheorem{lemma}{Lemma}[section]
\newtheorem{conj}{Conjecture}[section]
\begin{document}

\def\Z{\mathbb{Z}}
\def\R{\mathbb{R}}
\def\P{{\bf P}}
\def\E{{\bf E}}
\def\Var{{\bf Var}}
\def\cG{{\mathcal G}}
\def\cV{{\mathcal V}}
\def\cE{{\mathcal E}}
\def\coex{{\rm Coex}}
\def\main{{\rm main}}
\def\max{{\rm max}}
\def\Cox{\hfill \Box}
\def\eps{\varepsilon}

\title{Nonmonotonic coexistence regions for the \\
two-type Richardson model on graphs}

\author{Maria Deijfen \thanks{Department of Mathematics, Stockholm University,
106 91 Stockholm, Sweden. E-mail: mia@math.su.se.} \and Olle
H\"{a}ggstr\"{o}m
\thanks{Department of Mathematics, Chalmers University of Technology, 412 96 Gothenburg,
Sweden. E-mail: olleh@math.chalmers.se.}}

\date{February 2006}

\maketitle

\thispagestyle{empty}

\begin{abstract}

\noindent In the two-type Richardson model on a graph
$\cG=(\cV,\cE)$, each vertex is at a given time in state $0$, $1$
or $2$. A $0$ flips to a $1$ (resp.\ $2$) at rate $\lambda_1$
($\lambda_2$) times the number of neighboring $1$'s ($2$'s), while
$1$'s and $2$'s never flip. When $\cG$ is infinite, the main
question is whether, starting from a single $1$ and a single $2$,
with positive probability we will see  both types of infection
reach infinitely many sites. This has previously been studied on
the $d$-dimensional cubic lattice $\Z^d$, $d\geq 2$, where the
conjecture (on which a good deal of progress has been made) is
that such coexistence has positive probability if and only if
$\lambda_1=\lambda_2$. In the present paper examples are given of
other graphs where the set of points in the parameter space which
admit such coexistence has a more surprising form. In particular,
there exist graphs exhibiting coexistence at some value of
$\frac{\lambda_1}{\lambda_2} \neq 1$ and non-coexistence when this
ratio is brought closer to $1$.

\vspace{0.5cm}


\vspace{0.5cm}

\hspace*{5.95cm} Secondary 82B43.

\vspace{0.5cm}

\noindent Submitted March 1, 2005, final version accepted February
22, 2006.

\end{abstract}

\vfill\eject

\section{Introduction}  \label{sect:intro}

The two-type Richardson model on a graph $\cG=(\cV,\cE)$ is an
interacting particle system where at any time $t$, each vertex
$v\in \cV$ is in state $0$, $1$ or $2$. Here $1$ and $2$ may be
interpreted as two mutually exclusive types of infection, while
sites in state $0$ are thought of as not being infected. The
dynamics is that a site $v \in \cV$ in state $1$ or $2$ remains in
this state forever, while a $0$ flips to a $1$ (resp.\ $2$) at
rate $\lambda_1$ ($\lambda_2$) times the number of $1$'s ($2$'s)
among the neighbors of $v$, where two sites are said to be
neighbors if they share an edge in $\cE$. Here
$\lambda_1,\lambda_2>0$ are the two infection parameters of the
model. The graph $\cG$ will always be assumed to be countable and
connected.

This model has previously been studied on the $\Z^d$ lattice, that
is, on the graph whose vertex set is $\Z^d$ and whose edge set
consists of all pairs of sites at Euclidean distance $1$ from each
other. The main question is whether, when we start from a single
site in state $1$, a single site in state $2$, and all others
uninfected, we get positive probability for the event that both
types of infection succeed in reaching an infinite number of
sites. This event will in the following be referred to as {\em
infinite coexistence}. Note that, given the initial configuration,
the probability of infinite coexistence depends on $\lambda_1$ and
$\lambda_2$ only through their ratio
$\lambda=\frac{\lambda_2}{\lambda_1}$, as follows by a simple
time-scaling argument. For this reason, we may without loss of
generality set $\lambda_1=1$ and vary only $\lambda$
($=\lambda_2$). The following conjecture goes back to
H\"aggstr\"om and Pemantle \cite{HP1,HP2}.
\begin{conj}\label{conj:richardson_coex}
Infinite coexistence in the two-type Richardson model on $\Z^d$, $d\geq 2$,
starting from a single infected site of each type,
has positive probability if and only if $\lambda=1$.
\end{conj}
A good deal of progress has been made on this conjecture.
H\"aggstr\"om and Pemantle \cite{HP1} showed that for $d=2$ and
$\lambda=1$, infinite coexistence has positive probability. This
result was recently extended to $d \geq 3$ (as well as to more
general models) by Garet and Marchand \cite{GM} and independently
by Hoffman \cite{Hof}. As far as excluding infinite coexistence
for $\lambda\neq 1$, the best result to date is the following.
\begin{theorem}  \label{thm:HP2}
For the two-type Richardson model on $\Z^d$, $d\geq 2$, infinite
coexistence has probability $0$ for all but at most countably many
values of $\lambda$.
\end{theorem}
An analogous result for a related continuum model was obtained by
Deijfen et al.\ \cite{DHB}. Also, Deijfen and H\"aggstr\"om \cite{DH}
showed that the initial configuration, as long as there is a finite
nonzero number of infected sites of each type and one infection has not
already ``strangled'' the other, does not matter for the issue of
whether or not infinite coexistence has positive probability.
It is our experience from talking to colleagues about Theorem
\ref{thm:HP2} (and its continuum analogue) that they tend to react
with surprise at how weak this result is, and suggest that it
should be easy to improve in such a way as to obtain the ``only if''
direction of Conjecture \ref{conj:richardson_coex}. Their argument
is invariably the following.
\begin{equation}   \label{eq:monotonicity_intuition}
\begin{array}{l}
\mbox{Suppose for contradiction that infinite coexistence has} \\
\mbox{positive probability for some $\lambda > 1$. Then we can, due} \\
\mbox{to Theorem \ref{thm:HP2}, find some $\lambda' \in (1,
\lambda)$
for which infinite} \\
\mbox{coexistence has probability $0$. But this is absurd, since} \\
\mbox{surely it must be easier to get infinite coexistence if we} \\
\mbox{pick $\lambda$ closer to the symmetry point $1$.}
\end{array}
\end{equation}
Although we do agree with this intuition (see also Lundin \cite{L}
for some numerical evidence in support for it),
we think on the other hand that the claimed
monotonicity may not be easy to prove. In particular, we do not believe
that it is possible to establish using abstract arguments that
disregard the particular geometry of the $\Z^d$ lattice.

It is the purpose of this paper to support this point of view by
giving examples of other graphs, where the two-type Richardson
model behaves in a way that conflicts with the intuition about
monotonicity in $\lambda$. These graphs differ from the
$\mathbb{Z}^d$-lattice in that they are highly non-symmetric:
certain parts of the graph are designed specifically with
propagation of type-1 infection in mind, while other (different)
parts are meant for type-2 infection.

For a graph $\cG$, write $\coex(\cG)$ for the set of all $\lambda
\geq 1$ such that there exists an initial configuration $\xi \in
\{0,1,2\}^{\cV}$ which has only finitely many infected sites of
each type and for which the two-type Richardson model starting
from $\xi$ yields infinite coexistence with positive probability.
(Note that by time-scaling and interchange of $1$'s and $2$'s,
coexistence is possible for $\lambda$ if and only if it is
possible for $\lambda^{-1}$; hence no information is lost by
restricting to $\lambda \geq 1$.) In Sections \ref{sect:basic} and
\ref{sect:further}, we will exhibit examples of graphs $\cG$ that
demonstrate that, among others, the following kinds of coexistence
sets $\coex(\cG)$ are possible:
\begin{itemize}
\item
For any positive integer $k$, $\coex(\cG)$ may consist of exactly
$k$ points.
\item
$\coex(\cG)$ may be countably infinite.
\item
$\coex(\cG)$ may be an interval $(a,b)$ with $1<a<b$.
\end{itemize}
Note that all three examples show that the monotonicity intuition
suggested in (\ref{eq:monotonicity_intuition}) fails for general
graphs. However, as mentioned above, all examples will be highly
non-symmetric. A reasonable guess is that the intuition in
(\ref{eq:monotonicity_intuition}) is in fact correct on transitive
graphs.

As a complement to these examples, we will end the paper by giving
a positive result (Theorem \ref{thm:no_exotic_coex}) ruling out
a large class of more exotic coexistence regions, such
as those that are uncountable with zero Lebesgue measure.

Before moving on to the examples and results, let us say a few
words about the construction and the well-definedness of the
two-type Richardson model on graphs. To this end, let $S^1_t$ and
$S^2_t$ denote the set of type 1 and 2 infected vertices
respectively at time $t$, and, for a vertex set $A\in\cV$, write
$\partial A$ for the set of edges with one endpoint in $A$ and one
endpoint in $A^c$. One way to construct the model with parameters
$\lambda_1, \lambda_2>0$ is to assign i.i.d.\ exponential random
variables $\{X(e)\}_{e \in \cE}$ with mean 1 to the edges of $\cG$
and update the sets $S^1_t$ and $S^2_t$ inductively at discrete
time points $\{T_n\}$ as follows:

\begin{itemize}
\item[1.] Define $T_0=0$ and pick two bounded initial sets $S_0^1$
and $S_0^2$.

\item[2.] For $n\geq 1$, given $T_{n-1}$, $S_{T_{n-1}}^1$ and
$S_{T_{n-1}}^2$, define $T_n=\min\{T_n^1, T_n^2\}$, where
$$
T_n^i=\inf\{\lambda_iX(e);\,e\in\partial S_{T_{n-1}}^i\backslash
\partial S_{T_{n-1}}^j\},
$$
with $i,j\in\{1,2\}$ and $i\neq j$.

\item[3.] Let $S_t^i=S_{T_{n-1}}^i$ for $t\in[T_{n-1},T_n)$ and
$i=1,2$. Then, at time $T_n$, the sets are updated in that
infection is transferred through the edge defining $T_n$. More
precisely, if $T_n=T_n^1$ and $x$ is the uninfected end of the
edge where the infimum in the definition of $T_n^1$ is attained,
then $S_{T_n}^1=S_{T_{n-1}}^1\cup \{x\}$ and
$S_{T_n}^2=S_{T_{n-1}}^2$. Similarly, if $T_n=T_n^2$ then
$S_{T_{n-1}}^2$ is updated analogously, while $S_{T_{n-1}}^1$ is
left unchanged.

\end{itemize}

If $\cG$ has bounded degree (as all our examples will) and if
initially only finitely many vertices are infected, then it is
straightforward to see that almost surely no explosion will occur
(where explosion means that infinitely many transitions take place
in finite time), and that the process is Markovian with the
desired infection intensities.

\section{Basic examples}  \label{sect:basic}

We begin with our simplest example: a graph $\cG$ that admits
infinite coexistence in the two-type Richardson model if and only
if $\lambda$ is $\frac{1}{2}$ or $2$.
\begin{prop}  \label{prop:simplest_example}
There exists a graph $\cG$ with $\coex(\cG)= \{2\}$.
\end{prop}
{\bf Proof.} The graph we use to prove this result will look a bit
like a ladder: with two ``spines'' linked by a number of bridges.

More specifically, let $\{v_{1,j}\}_{j\geq 0}$ and
$\{v_{2,j}\}_{j\geq 0}$ be two sequences of vertices, each
internally linked by edges $\{e_{1,j}\}_{j \geq 0}:=\{\langle
v_{1,j}, v_{1,j+1}\rangle\}_{j \geq 0}$ and $\{e_{2,j}\}_{j \geq
0}:=\{\langle v_{2,j}, v_{2,j+1}\rangle\}_{j \geq 0}$,
respectively. These two infinite paths will be linked to each
other by finite paths, called bridges, where the $n$:th such
bridge emanates from $v_{1, a_n}$ and arrives at $v_{2, 2a_n}$.
Here $(a_1, a_2, \ldots)$ is a rapidly increasing sequence of
positive integers (how rapidly will be indicated later). The
$n$:th bridge will be called $B_n$ and have length $\left\lceil
a_n^{7/8} \right\rceil$, where $\lceil \cdot \rceil$ denotes
rounding up to the nearest integer.

Now consider the two-type Richardson model on this graph with
infection rates $\lambda_1=1$ and $\lambda_2=\lambda>0$ starting
with a single infected $1$ at $v_{1,0}$, and a single infected $2$
at $v_{2,0}$. It is easy to see that if a site on the first spine
$\{v_{1,j}\}_{j\geq 0}$ is ever infected by the type $2$
infection, then type $1$ is strangled (that is, it is cut off from
the possibility of ever infecting more than a finite number of
sites). Hence, the event $C_1$ of infinite growth of the type $1$
infection happens if and only if {\em all} sites on the first
spine are eventually infected by type $1$. Similarly, the event
$C_2$ of infinite growth of the type $2$ infection happens if and
only if all sites on the second spine are eventually infected by
type $2$.

With the edge representation indicated at the end of Section \ref{sect:intro},
let $D_{1,n}$ denote the event that
\begin{equation}  \label{eq:two_intercepts_one}
\sum_{j=0}^{a_n-1} X(e_{1,j}) \, > \,
\sum_{j=0}^{2a_n-1} \lambda^{-1}X(e_{2,j}) +
\sum_{e \in B_n} \lambda^{-1}X(e) \, .
\end{equation}
Note that, unless the type $1$ infection has already managed to
infect some site on the second spine before reaching $v_{1, a_n}$,
the site $v_{1, a_n}$ gets type $1$ infected if and only if
$D_{1,m}$ does {\em not} happen for any $m\leq n$. Analogously,
define $D_{2,n}$ as
\begin{equation}  \label{eq:one_intercepts_two}
\sum_{j=0}^{2a_n-1} \lambda^{-1}X(e_{2,j}) \, > \,
\sum_{j=0}^{a_n-1} X(e_{1,j}) + \sum_{e \in B_n} X(e) \, ,
\end{equation}
and note that, if type $2$ has not already infected some site on
the first spine before reaching $v_{2, 2a_n}$, then $v_{2, 2a_n}$
gets type $2$ infected if and only if $D_{2,n}$ does {\em not}
happen. Hence, the event $C= C_1 \cap C_2$ of infinite coexistence
happens if and only if none of the events $D_{1,1}, D_{1,2},
\ldots$ and $D_{2,1}, D_{2,2}, \ldots$ happen.

To get a grip on the probabilities of these events, it is useful to
introduce the variables
\[
T_{1,n} \, = \, \sum_{j=0}^{2a_n-1} \lambda^{-1}X(e_{2,j}) +
\sum_{e \in B_n} \lambda^{-1}X(e) - \sum_{j=0}^{a_n-1} X(e_{1,j})
\]
and
\[
T_{2,n} \, = \, \sum_{j=0}^{a_n-1} X(e_{1,j}) + \sum_{e \in B_n} X(e)
- \sum_{j=0}^{2a_n-1} \lambda^{-1}X(e_{2,j})
\]
and note that $D_{1,n}=\{T_{1,n} < 0\}$ and $D_{2,n}=\{T_{2,n} < 0\}$.
The expectation and variance of $T_{1,n}$ are
\[
\E\left[T_{1,n}\right] \, = \, (2\lambda^{-1}-1) a_n +
\lambda^{-1}\left\lceil a_n^{7/8} \right\rceil
\]
and
\[
\Var\left[T_{1,n}\right] \, = \, (2\lambda^{-2}+1) a_n +
\lambda^{-2}\left\lceil a_n^{7/8} \right\rceil \, ,
\]
and for $T_{2,n}$ we get
\[
\E\left[T_{2,n}\right] \, = \, (1-2\lambda^{-1})a_n +
\left\lceil a_n^{7/8} \right\rceil
\]
and
\[
\Var \left[T_{2,n}\right] \, = \,
(1+2\lambda^{-2})a_n + \left\lceil a_n^{7/8} \right\rceil \, .
\]

There are now three cases to consider separately, namely $\lambda>2$,
$\lambda=2$ and $\lambda<2$.

We begin with the case $\lambda>2$. We then have $2\lambda^{-1}-1
<0$, meaning that $\E\left[T_{1,n}\right]<0$ for $a_n$ large
enough. Writing $\neg$ for set complement, Chebyshev's inequality
gives
\begin{eqnarray}  \nonumber
\P\left(\neg D_{1,n} \right) & = & \P\left(T_{1,n} \geq 0\right) \\
\nonumber
& \leq & \frac{\Var\left[T_{1,n}\right]}{(\E\left[T_{1,n}\right])^2} \\
& = & \frac{(2\lambda^{-2}+1) a_n + \lambda^{-2}\left\lceil
a_n^{7/8} \right\rceil}{ \left((2\lambda^{-1}-1) a_n +
\lambda^{-1}\left\lceil a_n^{7/8} \right\rceil \right) ^2 },
\label{eq:first_Cheb_application}
\end{eqnarray}
which tends to $0$ as $a_n \rightarrow \infty$, and, since $a_n$
tends to $\infty$ with $n$, it follows that $\lim_{n \rightarrow
\infty}\P\left[D_{1,n} \right]=1$. Hence
\begin{eqnarray*}
\P(C) & \leq & \P(C_1) \\
& = & 1 - \P \left( \bigcup_{n=1}^\infty D_{1,n} \right) \\
& \leq & 1 - \lim_{n \rightarrow \infty}\P \left(D_{1,n} \right) \\
& = & 0 \, .
\end{eqnarray*}

In the case $\lambda<2$, we have $1-2\lambda^{-1}<0$ and hence
$\E\left[T_{2,n}\right]<0$ for $a_n$ large enough, whence
analogously to (\ref{eq:first_Cheb_application}) we get
\[
\P\left(\neg D_{2,n} \right) \, \leq \,
\frac{ (1+2\lambda^{-2})a_n + \left\lceil a_n^{7/8} \right\rceil }{ \left(
(1-2\lambda^{-1})a_n + \left\lceil a_n^{7/8} \right\rceil \right)^2} \, .
\]
This tends to $0$ as $a_n\rightarrow\infty$, so that $\lim_{n
\rightarrow\infty}\P(D_{2,n})=1$, implying that
\begin{eqnarray*}
\P(C) & \leq & \P(C_2) \\
& = & 1 - \P \left( \bigcup_{n=1}^\infty D_{2,n} \right) \\
& \leq & 1 - \lim_{n \rightarrow \infty}\P \left(D_{2,n} \right) \\
& = & 0 \, .
\end{eqnarray*}

The final case $\lambda=2$ is slightly more subtle. Both $\E[T_{1,n}]$ and
$\E[T_{2,n}]$ are then positive for any $n$, and another application of
Chebyshev's inequality gives
\begin{eqnarray*}
\P\left(D_{1,n} \right) & = & \P\left(T_{1,n} < 0\right) \\
\nonumber
& \leq & \frac{\Var\left[T_{1,n}\right]}{(\E\left[T_{1,n}\right])^2} \\
& = & \frac{(2^{-1}+1) a_n + 2^{-2}\left\lceil a_n^{7/8}
\right\rceil}{ \left( 2^{-1}\left\lceil a_n^{7/8} \right\rceil
\right) ^2 }
\end{eqnarray*}
which tends to $0$ as $a_n \rightarrow\infty$. Similarly,
$\P\left(D_{2,n} \right) \rightarrow 0$ as $a_n
\rightarrow\infty$. So far we have not specified how quickly the
numbers $a_n$ tends to $\infty$ with $n$. We are therefore free to
choose $(a_1, a_2, \ldots)$ in such a way that
\[
\sum_{n=1}^\infty \left( \P(D_{1,n}) + \P(D_{2,n}) \right) < 1 \, .
\]
With such a choice of $(a_1, a_2, \ldots)$, we get
\begin{eqnarray*}
\P(C) & = & 1- \P( \neg C_1 \cup \neg C_2 ) \\
& = & 1 - \P \left( \bigcup_{n=1}^\infty D_{1,n}
\cup \bigcup_{n=1}^\infty D_{2,n} \right) \\
& \geq & 1 - \sum_{n=1}^\infty \left( \P(D_{1,n}) + \P(D_{2,n}) \right) \\
& > & 0 \, .
\end{eqnarray*}
Having worked through the three cases $\lambda>2$, $\lambda=2$ and
$\lambda<2$, we have now shown that with the given initial
condition (infection $1$ at $v_{1,0}$ and infection $2$ at
$v_{2,0}$), infinite coexistence has positive probability if and
only if $\lambda=2$. In order to prove Proposition
\ref{prop:simplest_example}, it remains to show that no other
finite initial condition can yield infinite coexistence for any
other $\lambda \geq 1$.

Of course, if infections $1$ and $2$ switch places in the above
initial condition, then infinite coexistence has positive
probability if and only if $\lambda = \frac{1}{2}$. For other
initial conditions, note that infinite coexistence implies that
either
\begin{description}
\item{(i)} infection $1$ finds a path to $\infty$ which from some
point onwards belongs to the first spine, and infection $2$
similarly reaches $\infty$ along the second spine, or \item{(ii)}
vice versa.
\end{description}
But the above analysis shows that scenario (i) requires $\lambda=2$,
and similarly scenario (ii) requires $\lambda= \frac{1}{2}$, so
$\coex(\cG)= \{ 2 \}$ as desired. $\Cox$

\medskip\noindent
The reason we get a different behavior for the competition process
here as compared to on $\mathbb{Z}^d$ is related to the lack of
symmetry of the above graph. The two infinite spines provide two
separate paths to infinity on which the infection types can keep
in step with each other if their intensities are chosen to match
the density of vertices on the spines. Note also that, on the
above graph, the initial configuration is indeed important for the
possibility of infinite coexistence. For instance, as pointed out
by the end of the proof of Proposition
\ref{prop:simplest_example}, switching two single sources might
change the coexistence probability. Again, this is related to the
lack of symmetry of the graph and contrasts with the
$\mathbb{Z}^d$ case (on $\mathbb{Z}^d$, switching two sources
located at $x$ and $y$ respectively does not affect the
coexistence probability, since there exists an automorphism of the
graph $\mathbb{Z}^d$ that exchanges $x$ and $y$).

Of course, for any $\alpha \geq 1$, we can modify the
above construction by letting bridges connect $v_{1, a_n}$ and
$v_{2, \lceil \alpha a_n \rceil}$ rather than $v_{1, a_n}$ and
$v_{2, 2 a_n}$, thereby obtaining a graph $\cG$ with
$\coex(\cG)= \{ \alpha\}$. See Theorem \ref{thm:finite_coex_region}
for a more general result.

Next, we show how to turn $\coex(\cG)$ into an entire interval.
\begin{prop}  \label{prop:entire_interval}
There exists a graph $\cG$ for which $\coex(\cG)$ equals the
interval $[2,5]$.
\end{prop}
{\bf Proof.} As in the proof of Proposition
\ref{prop:simplest_example}, we take $\cG$ to consist of two
spines $\{v_{1,j}\}_{j\geq 0}$ and $\{v_{2,j}\}_{j\geq 0}$,
together with a sequence of bridges between them. This time, we
take the $n$'th bridge $B_n$ to begin at $v_{1, a_n}$ and end at
$v_{2, 4a_n}$, and to have length $a_n + \left\lceil a_n^{7/8}
\right\rceil$.

Note that (i) and (ii) in the proof of Proposition
\ref{prop:simplest_example} are the only two possible scenarios
for infinite coexistence, and that it is therefore sufficient to
consider the initial condition with a single $1$ at $v_{1,0}$ and
a single $2$ at $v_{2,0}$. With this initial condition, the same
kind of applications of Chebyshev's inequality as in the proof of
Proposition \ref{prop:simplest_example} show that
$\P(\cap_{n=1}^\infty D_{1,n}^c)= 0$ for $\lambda>5$, that
$\P(\cap_{n=1}^\infty D_{2,n}^c)= 0$ for $\lambda<2$, and that
$\P(C)>0$ for $\lambda \in [2,5]$ provided the sequence $(a_1,
a_2, \ldots)$ grows sufficiently fast. Thus, $\coex(\cG)= [2,5]$.
$\Cox$

\medskip\noindent
Several variations of Proposition \ref{prop:entire_interval} are
easily obtained, such as the following
(and combinations thereof):
\begin{enumerate}
\item For any $1\leq \alpha < \alpha'$, the coexistence region
$[2,5]$ can be replaced by $\coex(\cG)= [\alpha,\alpha']$ by
noting that, if the bridge $B_n$ is taken to begin at $v_{1,a_n}$
and end at $v_{2,\left\lceil\gamma a_n\right\rceil}$ and to have
length $\left\lceil\beta a_n+ a_n^{7/8}\right\rceil$, then infinite
coexistence is possible for
$\lambda\in[\gamma/(1+\beta),\gamma+\beta]$. The
 interval $[\alpha,\alpha']$ is hence obtained by letting $B_n$ end at
\[
v_{2,
\left\lceil\frac{\alpha(1+\alpha')a_n}{1+\alpha}\right\rceil} \, ,
\]
and have length
$\left\lceil\frac{(\alpha'-\alpha)a_n}{1+\alpha}\right\rceil +
\left\lceil a_n^{7/8} \right\rceil$. \item If the $B_n$'s are
taken to have length $a_n - \left\lceil a_n^{7/8} \right\rceil$
rather than $a_n + \left\lceil a_n^{7/8} \right\rceil$, then the
coexistence region $[2, 5]$ is replaced by the open interval
$(2,5)$. \item If the lengths of the $B_n$'s are taken to be $a_n$
(with no lower order correction), start from $v_{1, a_n}$ and end
at
\[
v_{2, 3a_n+ \left\lceil a_n^{7/8} \right\rceil} \, ,
\]
then we get $\coex(\cG) = (2, 5]$. The other half-open interval
$[2, 5)$ can be obtained analogously.
\end{enumerate}

\section{Further examples}  \label{sect:further}

The following result generalizes Proposition \ref{prop:simplest_example}.
\begin{theorem}  \label{thm:finite_coex_region}
For any $k$ and any $\alpha_1, \ldots ,\alpha_k \in [1, \infty)$,
there exists a graph $\cG$ with $\coex(\cG)=\{\alpha_1, \ldots,
\alpha_k\}$.
\end{theorem}
Intuitively, it is not hard to figure out what kind of
modification of the example in Proposition
\ref{prop:simplest_example} that would lead to Theorem
\ref{thm:finite_coex_region}. Let $\cG$ contain $k+1$ spines,
whose respective vertex sets we may denote $\{v_{\main,j}\}_{j\geq
0}$ and $\{v_{1,j}\}_{j\geq 0}, \ldots, \{v_{k,j}\}_{j\geq 0}$. As
before, $(a_1, a_2, \ldots)$ will be a rapidly increasing
sequence, and for each $n$ there will be bridges $B_{1,n}, \ldots,
B_{k,n}$, each with length $\left\lceil a_n^{7/8} \right\rceil$,
bridge $B_{m,n}$ starting at $v_{\main, a_n}$ and ending at $v_{m,
\left\lceil \alpha_m a_n \right\rceil }$.

Infinite coexistence should now be possible for $\lambda=\alpha_m$
by means of infection $1$ taking over the main spine
$\{v_{\main,j}\}_{j\geq 0}$, and infection $2$ taking over the
$m$'th auxiliary spine $\{v_{m,j}\}_{j\geq 0}$ (while all other
spines are conquered by infection $1$). On the other hand, no such
scenario seems possible when $\lambda \geq 1$ is not an element of
$\{\alpha_1, \ldots ,\alpha_k\}$.

Proving this turns out to be a technically somewhat more
challenging task compared to what we did in Section
\ref{sect:basic}, the reason being that there is no slick
decription of the possible ways of infinite coexistence like (i)
and (ii) in the proof of Proposition \ref{prop:simplest_example}.
On the contrary, for $k\geq 2$ we may (in principle) imagine an
infinite coexistence scenario where each infection ``zig-zags''
between the main spine and the other spines in a more or less
complicated manner. For this reason, we choose to construct the
graph in a more iterative fashion, and to base our arguments on
the following lemma. For the two-type Richardson model on a graph
$\cG=(\cV, \cE)$ and a vertex $v\in \cV$, write $T_v$ for the time
at which $v$ becomes infected.
\begin{lemma}  \label{lem:for_the_technical_constr}
Consider the two-type Richardson model with infection rates
$\lambda_1=1$ and $\lambda_2 = \lambda \geq 1$ on a graph $\cG$
constructed as follows: Let $\cG'=(\cV', \cE')$ be an arbitrary
finite graph with $k \geq 1$ distinguished vertices $v_1, \ldots,
v_k \in \cV'$, and obtain $\cG$ by, for $i=1, \ldots, k$,
attaching to $v_i$ an infinite path with vertex set
$\{v_{i,j}\}_{j\geq 1}$.

Fix $\eps>0$. Then, for all sufficiently large $N$ (depending on
$\eps$ but not on $\lambda$),
we have that for any $\xi \in \{1,2\}^{\{v_1, \ldots, v_k\}}$ and any
initial condition such that
\begin{description}
\item{(a)} only vertices in $\cG'$ are initially infected, and
\item{(b)} the event $A_\xi$ that the infection proceeds in such a way
that for each $i \in \{1, \ldots, k\}$, vertex $v_i$ eventually
gets infection $\xi(v_i)$, has positive probability,
\end{description}
the following holds for each $i= \{1, \ldots, k\}$:
\[
\P \left(
T_{v_{i, N}} \in
\left( \lambda^{1 - \xi(v_i)} N- N^{3/4},
\lambda^{1 - \xi(v_i)}N+N^{3/4} \right) \, \Big| \, A_\xi \right)
\, \geq \, 1 - \eps \, .
\]
\end{lemma}
{\bf Proof.} If the initial condition is such that $v_1, \ldots,
v_k$ are all initially infected, then the lemma is immediate from
the Central Limit Theorem applied to the sum
$$
X(\langle v_i, v_{i,1} \rangle) + \sum_{j=1}^{N-1} X(\langle
v_{i,j}, v_{i,j+1} \rangle)
$$
for each $i \in \{1, \ldots, k\}$. For other initial conditions,
note that, since $\mathcal{G}'$ is finite and connected, the time
$T_{v_i}$ when the vertex $v_i$ is infected is almost surely
finite, implying that

\begin{equation}\label{eq:1}
P\left(T_{v_i}\leq N^{1/8}\Big|A_{\xi}\right)\geq 1-\eps/2
\end{equation}

\noindent for large $N$. On $A_{\xi}$, the vertex $v_i$ is
infected by the type $\xi(v_i)$ infection and the Central Limit
Theorem applied to the same sum as above hence gives that

\begin{equation}\label{eq:2}
\P \left( T_{v_{i, N}}-T_{v_i} \in \left( \lambda^{1 - \xi(v_i)}
N- N^{3/4}, \lambda^{1 - \xi(v_i)}N+N^{5/8} \right) \, \Big| \,
A_\xi \right) \, \geq \, 1 - \eps/2 \, .
\end{equation}

\noindent The desired estimate now follows by taking $N$ large
enough to ensure that both (\ref{eq:1}) and (\ref{eq:2}) hold for
all $i\in\{1,\ldots,k\}$.$\Cox$

\medskip\noindent
{\bf Proof of Theorem \ref{thm:finite_coex_region}.}
Fix the coexistence set $\{\alpha_1, \ldots \alpha_k\}$ that we are trying to
obtain, and define $\alpha_0 = 1$
and $\alpha_\max= \max\{\alpha_1, \ldots, \alpha_k\}$.
The graph $\cG=(\cV, \cE)$ that will
serve as an example, will be obtained from a sequence of graphs
$\{\cG_n=(\cV_n, \cE_n)\}_{n \geq 0}$ which is increasing in the sense
that
\[
\cV_0 \subseteq \cV_1 \subseteq \cdots
\]
and
\[
\cE_0 \subseteq \cE_1 \subseteq \cdots \, .
\]
The ''limiting'' graph $\cG=(\cV, \cE)$ is then given by $\cV=
\bigcup_{n=0}^\infty \cV_n$ and $\cE= \bigcup_{n=0}^\infty \cE_n$.
The graphs $\cG_n$ are defined inductively as follows.

Let $(b_1, b_2, \ldots)$ be a rapidly growing sequence of positive
integers; how rapidly $b_n \rightarrow \infty$ with $n$ will be
specified later. The construction begins with taking $\cG_0$ to be
the complete graph on $k+1$ vertices $x_{0,0}, x_{0,1}, \ldots,
x_{0, k}$. For the induction step, suppose that we have the graph
$\cG_n$ with $k+1$ distinguished vertices $x_{n,0},, x_{n,1},
\ldots, x_{n, k}$, and let $\cG_{n+1}$ be the graph obtained from
$\cG_n$ by the following amendments:
\begin{description}
\item{(a)}
for $i=0, \ldots, k$, attach to $x_{n,i}$ a path of length
$\left\lceil \alpha_i b_{n+1} \right\rceil$, denoting the last vertex
of this path by $x_{n+1,i}$, and
\item{(b)}
for $i=1, \ldots, k$,  link $x_{n+1, 0}$ and $x_{n+1, i}$ by a path
(called a bridge and denoted $B_{n+1, i}$) of length
$\left\lceil b_{n+1}^{7/8} \right\rceil$.
\end{description}
This defines $(\cG_0, \cG_1, \ldots)$ and $\cG$, apart from that
the sequence $(b_1, b_2, \ldots)$ has not been specified. In order
to make that choice, begin by specifying two decreasing sequences
$(\delta_1, \delta_2, \ldots)$ and $(\eps_1, \eps_2, \ldots)$ of
positive numbers tending to $0$, with the second sequence having
the additional property that
\begin{equation}  \label{eq:summability}
\sum_{n=1}^\infty \eps_n \, < \, \frac{1}{k+2} \, .
\end{equation}
Given $b_1, \ldots, b_n$ (and, thus, $\cG_1, \ldots, \cG_n$), pick
$b_{n+1}$ large enough so that the following conditions hold:
\begin{description}
\item{(i)} For $\xi \in \{1,2\}^{\{x_{n,0}, \ldots, x_{n, k}\}}$,
let $A_{n, \xi}$ denote the event that $x_{n,0}, \ldots, x_{n, k}$
are infected by type $\xi(x_{n,0}), \ldots, \xi(x_{n, k})$,
respectively. For an arbitrary such $\xi$ and arbitrary initial
conditions that are confined to $\cG_n$ and that makes it possible
for $A_{n,\xi}$ to happen before any vertices outside $\cG_n$ are
infected, we have for $i=0, \ldots, k$ and any $\lambda \geq 1$
that
\begin{eqnarray*}
\lefteqn{\hspace{-27mm} \P\left( T_{x_{n+1},i} \in
\left(\lambda^{1-\xi(x_{n,i})}\alpha_i b_{n+1} -\alpha_i
b_{n+1}^{3/4}, \lambda^{1-\xi(x_{n,i})} \alpha_i b_{n+1}+\alpha_i
b_{n+1}^{3/4} \right) \Big| A_{n,\xi} \right)}
\hspace{40mm}  \\
& \geq & 1 - \eps_{n+1}
\end{eqnarray*}
(note that this holds for $b_{n+1}$ large enough,
due to Lemma \ref{lem:for_the_technical_constr}).
\item{(ii)}
\[
2 \alpha_\max b_{n+1}^{3/4} < \frac{\delta_{n+1}b_{n+1}}{2}
\]
\item{(iii)}
\[
\mbox{ } \hspace{-6mm}
\P\left( \sum_{e\in B_{n+1,i}} X(e) \in
\left(2\alpha_\max b_{n+1}^{3/4},
\frac{\delta_{n+1}b_{n+1}}{2} \right) \mbox{ for } i=1, \ldots, k \right)
\, \geq \, 1 - \eps_{n+1}
\]
(note that this holds for $b_{n+1}$ large enough, since by the
Weak Law of Large Numbers the sum $\sum_{e\in B_{n+1,i}} X(e)$ is
concentrated around $b_{n+1}^{7/8}$ for large values of
$b_{n+1}$).
\end{description}
This specifies $\cG$. Now fix $i \in \{1, \ldots, k\}$, and
consider the two-type Richardson model with $\lambda_1=1$ and
$\lambda=\lambda_2 = \alpha_i$, starting with vertex $x_{0,i}$ in
state $2$, vertices $x_{0,0}, x_{0,1}, \ldots, x_{0, i-1}, x_{0,
i+1}, \ldots, x_{0,k}$ in state $1$, and all other vertices
uninfected.

For $n=1,2, \ldots$, let $D_{n,i}$ denote the event that $x_{n,0}$
eventually gets type 1 infected, that $x_{n,i}$ eventually gets
type 2 infected and that, at the time when the last one (in time)
of $x_{n,0}$ and $x_{n,i}$ is infected, there is no type 2
infection on any of the spines $j\neq i$. We have that
\[
\P(D_{1,i}) \, \geq \, 1 - (k+2) \eps_1 \, ,
\]
because the choice of $b_1$ implies that with probability at least
$1 - (k+2) \eps_1$, the type 1 infection reaches $x_{n,0}$ before
the type 2 infection has crossed the bridge $B_{1,i}$ and the type
2 infection reaches $x_{n,i}$ before the type 1 infection has
crossed $B_{1,i}$.

More generally, for any $n$, it follows from the choice of
$b_{n+1}$ that
\[
\P(D_{n+1, i} \, | \, D_{n,i}) \, \geq \, 1- (k+2)\eps_{n+1} \, .
\]
Hence, using (\ref{eq:summability}), we get
\begin{eqnarray*}
\P(C) & \geq & \P( \liminf_{n \rightarrow \infty} D_{n,i}) \\
& \geq &  \liminf_{n \rightarrow \infty} \P(D_{n,i}) \\
& \geq & \prod_{n=1}^\infty (1 - (k+2) \eps_n) \\
& \geq & 1 - \sum_{n=1}^\infty (k+2) \eps_n \\
& > & 0.
\end{eqnarray*}
We have thus shown that $\alpha_i\in \coex(\cG)$ for each $i \in
\{1, \ldots, k\}$.

What remains is to show that for any $\lambda \geq 1$ such that
$\lambda \not\in \{ \alpha_1, \ldots, \alpha_k\}$, and any finite
initial condition, we have that $\P(C)=0$. We will assume for the
moment that $\lambda>1$; fix such a $\lambda\not\in \{ \alpha_1,
\ldots, \alpha_k\}$ and an arbitrary initial condition, and pick
$n$ large enough so that
\begin{description}
\item{(a)} the
set of infected sites in the initial condition is confined to $\cG_n$
\item{(b)} $|\lambda - \alpha_i|> \delta_{n+1}$ for all
$i \in \{0, \ldots, k\}$, and
\item{(c)} $\frac{| \lambda - \alpha_i | - \delta_{n+1}}{\lambda} b_{n+1}
> (\alpha_i +1) b_{n+1}^{3/4}$  for all
$i \in \{0, \ldots, k\}$.
\end{description}
We now claim that
\begin{equation}  \label{eq:easy_to_get_killed}
\P(\mbox{both types of infection reach }\{x_{n+1,0}, \ldots, x_{n+1, k}\})
\, \leq \, (k+2)\eps_{n+1} \, .
\end{equation}
Once this is established we are done, because $\P(C)$ is bounded
by the left-hand side of (\ref{eq:easy_to_get_killed}) while the
right-hand side can be made arbitrarily small by picking $n$ even
larger.

To show that (\ref{eq:easy_to_get_killed}) holds, it suffices to show that
it holds even if we condition on the types of infection
$\xi= (\xi(x_{n,0}), \ldots, \xi(x_{n,k}))$ that reach
$x_{n,0}, \ldots, x_{n,k}$. Given such a $\xi$ (with both types of
infection appearing), define
\[
\alpha_\xi \, = \, \min \{ \alpha_i: \, i \in \{0,1, \ldots, k\},
\xi(x_{n,i})=2 \}
\]
and then fix $i$ in such a way that $\alpha_\xi= \alpha_i$. By the
choice of $n$ we have that $|\lambda- \alpha_\xi| > \delta_{n+1}$.
There are now two cases to consider: $\lambda> \alpha_\xi$ and
$\lambda< \alpha_\xi$. In the former case $\lambda> \alpha_\xi$,
we have by the choice of $b_{n+1}$ that with conditional
probability at least $1-(k+2)\eps_{n+1}$, infection $2$ claims all
of $x_{n+1,0}, \ldots, x_{n+1,k}$ (by rushing from $x_{n,i}$ to
$x_{n+1,i}$, then across the bridge $B_{n+1,i}$, and then across
the $k-1$ other bridges emanating from $x_{n+1, 0}$, all of this
before infection $1$ arrives at any of the $x_{n+1, j}$'s).
Similarly, in the latter case $\lambda< \alpha_\xi$, we get by the
choice of $b_{n+1}$ that with conditional probability at least
$1-(k+2)\eps_{n+1}$, infection $1$ claims all of $x_{n+1,0},
\ldots, x_{n+1,k}$ (by rushing from $x_{n,0}$ to $x_{n+1, 0}$ and
then across all bridges emanating from $x_{n+1, 0}$). Thus,
(\ref{eq:easy_to_get_killed}) is established when $\lambda>1$.

It only remains to deal with the case
$\lambda=1\not\in \{ \alpha_1, \ldots, \alpha_k\}$. This follows by a similar
argument: given $\xi = (\xi(x_{n,0}), \ldots, \xi(x_{n,k}))$,
the infection that has captured $x_{n,0}$ will,
with probability at least $1 - (k+2) \eps_{n+1}$, capture all of
$x_{n+1, 0}, \ldots, x_{n+1,k}$.
$\Cox$

\medskip\noindent
It now only takes a minor extension of the above construction to get
a graph that can be used to establish the following result.
\begin{theorem}  \label{thm:countably_infinite}
There exists a graph $\cG$ for which $\coex(\cG)$ is countably
infinite.
\end{theorem}
{\bf Proof.}
Take $(\alpha_1, \alpha_2, \ldots)$ to be an unbounded and
strictly increasing sequence with $\alpha_1\geq 1$. We will construct
a graph with $\coex(\cG)= \{\alpha_1, \alpha_2, \ldots\}$. As in the proof
of Theorem \ref{thm:finite_coex_region}, take $(b_1, b_2, \ldots)$ to be
a rapidly growing sequence of positive integers, and construct
$\cG= ( \cV, \cE)$ as a limit of an increasing sequence
$\{\cG_n = (\cE_n, \cV_n)\}_{n \geq 0}$ of finite graphs.

Take $\cG_0$ to consist of a single vertex $x_{0,0}$, and proceed
inductively: Given $\cG_n$, with $n+1$ distinguished vertices
$x_{n,0}, x_{n,1}, \ldots, x_{n,n}$, obtain $\cG_{n+1}$ by
decorating $\cG_n$ as follows.
\begin{description}
\item{(a)}
For $i=0,1,\ldots, n$, attach to $x_{n,i}$ a path of length
$\left\lceil \alpha_i b_{n+1} \right\rceil$ and denote the last
vertex of this path by $x_{n+1, i}$,
\item{(b)}
to the vertex $x_{n,n}$, attach an additional path of length
$\left\lceil \alpha_{n+1}b_{n+1} \right\rceil$, and denote the last
vertex of this path by $x_{n+1,n+1}$,
\item{(c)}
for $i=1, \ldots, n$, link $x_{n+1,0}$ and $x_{n+1, i}$ by a path of length
$\left\lceil b_{n+1}^{7/8} \right\rceil$.
\end{description}
To show that if $b_n \rightarrow \infty$ sufficiently fast as
$n\rightarrow \infty$, the graph $\cG$ gets coexistence region
$\{\alpha_1, \alpha_2, \ldots\}$, is now a completely
straightforward modification of the proof of Theorem
\ref{thm:finite_coex_region}. $\Cox$

\section{A positive result}  \label{sect:positive}

In the proof of Theorem \ref{thm:countably_infinite}, we required
that the candidate coexistence set $\{\alpha_1, \alpha_2,
\ldots\}$ could be written as an increasing unbounded sequence,
which is tantamount to saying that it has no accumulation points.
This condition is certainly not necessary for a countable set to
arise as a coexistence region for some graph, but we do not know
whether it can simply be removed.

One could ask for graphs with more exotic coexistence regions,
such as for instance examples whose coexistence regions are
uncountable with zero Lebesgue measure. That kind of behavior is,
however, ruled out by the following result. For a set $A \subseteq
\R$, write $A^c$ for its complement and $\partial A$ for its
boundary.
\begin{theorem}  \label{thm:no_exotic_coex}
For any graph $\cG$, the coexistence region $\coex(\cG)$ contains
at most countably many points in its boundary $\partial\coex(\cG)$.
\end{theorem}
{\bf Proof.} For a given finite initial condition $\xi$, write
$\coex_\xi(\cG)$ for the set of $\lambda$'s such that infinite
coexistence has positive probability with initial condition $\xi$.
Since $\cG$ is countable, there are only countably many finite
initial conditions $\xi$, and in order to prove the theorem it
therefore suffices to show for any $\xi$ that $\coex_\xi(\cG) \cap
\partial\coex_\xi(\cG)$ is countable.

Fix $\xi$. For the two-type Richardson model on $\cG$ with
parameters $\lambda_1=1$ and $\lambda_2=\lambda$, and initial
condition $\xi$, let $\theta_1(\lambda)$ denote the probability
that infection $1$ reaches only finitely many sites, and let
$\theta_2(\lambda)$ denote the probability that infection $2$
reaches {\em infinitely} many. We may assume that $\cG$ is
infinite (otherwise the statement of the theorem is trivial), in
which case we have $\theta_1 (\lambda) \leq \theta_2(\lambda)$,
with
\[
\left\{
\begin{array}{ll}
\theta_1(\lambda) = \theta_2 (\lambda)
& \mbox{if } \lambda \not\in \coex_\xi(\cG) \\
\theta_1(\lambda) < \theta_2 (\lambda)
& \mbox{if } \lambda \in \coex_\xi(\cG)  \, .
\end{array}  \right.
\]
The standard monotone coupling of the two-type Richardson model for different
$\lambda$'s (see, e.g., \cite{HP2}) shows that $\theta_1(\lambda)$
and $\theta_2(\lambda)$ are increasing functions of $\lambda$, whence
they can have at most countably many discontinuities.

Now define sets $A_-, A_+ \subseteq \R$ as
\[
A_- = \{\lambda \in \coex_\xi(\cG):
\, ( \lambda - \eps, \lambda) \not\subseteq
\coex_\xi(\cG) \mbox{ for all } \eps >0 \}
\]
and
\[
A_+ = \{\lambda \in \coex_\xi(\cG):
\, ( \lambda, \lambda + \eps ) \not\subseteq
\coex_\xi(\cG) \mbox{ for all } \eps >0 \} \, .
\]
For $\lambda^* \in A_-$, we see that
\[
\lim_{\lambda \nearrow \lambda^*} \theta_1 ( \lambda) \, = \,
\lim_{\lambda \nearrow \lambda^*} \theta_2 ( \lambda)
\]
by considering a subsequence of $\lambda$'s beloning to
$\neg\coex_\xi(\cG)$. Hence,
\begin{eqnarray*}
\lim_{\lambda \nearrow \lambda^*} \theta_2 ( \lambda)
& = & \lim_{\lambda \nearrow \lambda^*} \theta_1 ( \lambda) \\
& \leq &  \theta_1 ( \lambda^*) \\
& < & \theta_2 ( \lambda^*) \, ,
\end{eqnarray*}
so that $\theta_2(\lambda)$ has a discontinuity at $\lambda=\lambda^*$, and
it follows that $A_-$ is countable.

We similarly get for any $\lambda^* \in A_+$ that
$\theta_1(\lambda)$ has a discontinuity at $\lambda=\lambda^*$, so that
$A_+$ is countable as well. Note finally that
\[
\coex_\xi ( \cG) \cap \partial \coex_\xi ( \cG) \, = \,
A_- \cup A_+ \, ,
\]
which gives the desired conclusion that
$\coex_\xi ( \cG) \cap \partial \coex_\xi ( \cG)$ is countable.
$\Cox$

\end{document}